\documentclass[12pt]{amsart}
\usepackage[mathscr]{euscript}
\usepackage{amscd}
\theoremstyle{plain}
\newtheorem{theorem}{Theorem}[section]
\newtheorem{lemma}[theorem]{Lemma}
\newtheorem{proposition}[theorem]{Proposition}

\newtheorem*{EOWT}{Edge of the Wedge Theorem}

\theoremstyle{definition}

\newtheorem{definition}[theorem]{Definition}

\theoremstyle{remark}
\newtheorem{remark}[theorem]{Remark}

\numberwithin{equation}{section}

\newcommand{\V}{\mathcal V}
\newcommand{\D}{\mathcal D}
\newcommand{\N}{\mathcal N}
\newcommand{\M}{\mathcal M}
\newcommand{\F}{\mathbb F}
\newcommand{\C}{\mathbb C}
\newcommand{\R}{\mathbb R}
\newcommand{\bP}{\mathbb P}
\newcommand{\zb}{\overline z}
\newcommand{\mapsdownto}{\begin{array}c\scriptstyle\top\\[-8pt]
                                       \Big\downarrow\end{array}}
\begin{document}

\title[Involutive Structure]
{The Involutive Structure on the\\[8pt] Blow-Up of $\R^n$ in $\C^n$}
\author[M. Eastwood]{Michael Eastwood}
\address{
  Department of Pure Mathematics\\
  University of Adelaide\\
  South AUSTRALIA 5005}
\email{meastwoo@spam.maths.adelaide.edu.au}

\author[R. Graham]{C. Robin Graham}
\address{
  Department of Mathematics\\
  University of Washington \\
  Box 354350\\
  Seattle, WA 98195-4350}
\email{robin@math.washington.edu}
\maketitle

\thispagestyle{empty}

\renewcommand{\thefootnote}{}
\footnotetext{This research was supported by the Australian Research
Council and Department of Science and Technology and NSF grant
\#DMS-9303497.  This support and the hospitality extended
to the first author by the University of Washington are gratefully
acknowledged.}

\section{Introduction}\label{intro}
Many interesting geometric structures can be defined by specifying a
smooth subbundle $\V$ of the complexified tangent bundle of the underlying
smooth manifold, subject to an integrability
condition.  Examples include foliations, complex structures, and CR
structures.  In general such structures are called involutive,
or formally integrable, and their study has been the starting point of
far-reaching general investigations
(see, for example, \cite{t}, \cite{ct}, and \cite{hj}).  However, most
naturally occurring examples, for instance all those mentioned above, have the
property that $\V \cap \overline{\V}$ has constant rank.  In recent work on
integral geometry (\cite{begm}, \cite{e}, and \cite{be}), natural examples of
involutive structures have arisen for which the rank of $\V \cap
\overline{\V}$ changes along a hypersurface.  For these examples the
underlying manifold is the real blow-up of $\R\bP^n$ in $\C\bP^n$ for
various $n$. In this article we consider these new involutive structures
from an analytic point of view.

Locally we may as well consider the blow-up of $\R^n$ in $\C^n$.
This blow-up $B$ is a smooth real $2n$-manifold with a distinguished
hypersurface $\Sigma$,
the inverse image of $\R^n$ under the blow-down map
$b: B \rightarrow \C^n$.  The blow-up is defined
precisely so that the image under $b$ of a neighborhood of a point
of $\Sigma$ is a localized wedge in $\C^n$, i.e.\ the product of an open
set in $\R^n$ with a localized cone in $i\R^n$.
Away from $\Sigma$, $b$ is a diffeomorphism and the involutive structure is
just the lift of the complex structure on $\C^n \setminus \R^n$;
the bundle $\V$ of $(0,1)$ vectors extends smoothly across $\Sigma$
but $\dim(\V \cap \overline{\V}) = n-1$ there.

A solution of an involutive structure is a function or
distribution annihilated by all sections of $\V$.  An important analytic
problem is to understand the regularity of solutions of a given
involutive structure.  We show that
the involutive structure on the blow-up of $\R^n$ in $\C^n$ is hypocomplex,
which means that any solution is locally a holomorphic function of a basic set
of independent solutions.  {From} the point of view of the theory of involutive
structures, hypocomplexity is the strongest possible regularity property
for solutions; complex structures are hypocomplex.  Our first proof of the
hypocomplexity is an elementary power series argument.  We also give
another argument using the Edge of the Wedge Theorem of several complex
variables.  This is quite straightforward: a solution of the involutive
structure near a point of $\Sigma$ defines a holomorphic function on a
localized wedge and the Edge of the Wedge Theorem provides the extension
needed to deduce hypocomplexity.  It turns out that this argument can be
reversed; the Edge of the Wedge Theorem is an easy consequence of the
hypocomplexity of the involutive structure.  Thus, our power series argument
provides a relatively simple new proof of the Edge of the Wedge Theorem.
The idea of introducing the blow-up in this context seems to us
particularly natural; the wedge blows up to an open set with boundary
values along a hyperplane, where standard arguments can be applied.

In \S\ref{invstr} we describe in detail the involutive structure
on the blow-up of $\R^n$ in $\C^n$ and give the power series argument for
hypocomplexity.  We also use a similar power series argument to show
that the associated inhomogeneous equations have infinite dimensional first
local cohomology by showing that inhomogeneous terms of a particular form must
be real-analytic if there is a solution.  In \S\ref{eowt} we explain the
relationship with the Edge of the Wedge Theorem.  Finally in
\S\ref{intgeom} we review the integro-geometric considerations from
\cite{begm}, \cite{e}, and \cite{be} which
gave rise to the compact version of the involutive structure on the blow-up
in the first place.

\section{The Involutive Structure}\label{invstr}
We begin by reviewing some of the basic notions which will be relevant for
us concerning involutive
structures; see \cite{t} for
elaboration.  An involutive structure on a smooth
manifold $M$ is a smooth complex subbundle $\V \subset \C TM$
satisfying the formal integrability condition $[\V,\V]
\subset \V$.  Set $d=\dim_{\R} M - \dim_{\C} \V$.
A {\em solution} of $\V$ is a distribution $f \in \D'(M)$
satisfying $Lf=0$ for all sections $L$ of $\V$.
Familiar special cases include the following:
\begin{enumerate}
\item $\V = \overline{\V}$.  In this case $\V$ is the complexification of a
real Frobenius-integrable distribution on $M$, so $\V$ defines a foliation.
Locally one can find $d$ independent smooth real-valued solutions $x^1,
\ldots, x^d$
of~$\V$; any solution of $\V$ is of the form $g(x^1, \ldots, x^d)$ for
some \mbox{$g \in \D'(\R^d)$}.  Such involutive structures are called {\em
real}.
\item $\V \cap \overline{\V} = \{0\}$ and $\dim \V = \frac{1}{2} \dim M$.  In
this case $\V$ defines an integrable almost complex structure on $M$.  By
the Newlander-Nirenberg Theorem, locally one can find $d$
independent smooth solutions $z^1, \ldots, z^d$
of $\V$; then any solution of $\V$ is of the form
$h(z^1, \ldots, z^d)$ for
some holomorphic function~$h$.  Such involutive structures are called
{\em complex}.
\end{enumerate}

The involutive structure $\V$ is said to be {\em locally integrable} at $m \in
M$ if in some neighborhood of $m$ one can find $d$ smooth solutions
$f^1, \ldots, f^d$ of $\V$ with $df^1, \ldots, df^d$ linearly independent.
$\V$ is said to be locally integrable if it is locally integrable at every
$m \in M$.  The examples above are clearly locally integrable.  It
follows from the holomorphic version of the Frobenius
Theorem by complexification that any real-analytic involutive structure on
a real-analytic manifold is locally integrable.

An involutive structure is said to be a (generalized)   CR structure if
$\V \cap \overline{\V} = \{0\}$.
If additionally $\dim \V = \frac{1}{2}(\dim M- 1)$,
$\V$ is said to be CR of hypersurface type.  If such a $\V$ is locally
integrable at $m \in M$, then the image of the map
$f=(f^1, \ldots, f^d): M \to \C^d$ is a smooth hypersurface of $\C^d$ and
$f_*(\V)=T^{0,1}\C^d \cap \C Tf(M)$; one says that the CR structure is
embeddable.  There are CR structures which are not embeddable; these
provide examples of involutive structures which are not locally
integrable.

One of our primary interests is the regularity of solutions of locally
integrable involutive structures.
Complex involutive structures exhibit the
best possible regularity behavior:  every distribution solution is a
holomorphic function of a basic set of solutions.  The hypocomplex
involutive structures are defined by this regularity property:

\begin{definition}
Let $\V$ be an involutive structure on $M$, locally integrable at $m \in M$,
with $f^1, \ldots, f^d$ a set of independent smooth solutions of $\V$ on a
fixed neighborhood of $m$ as above.  $\V$ is said to be hypocomplex at $m$
if every solution of $\V$ defined in some neighborhood of $m$ is (on a possibly
smaller neighborhood) of the form
$h(f^1, \ldots, f^d)$ for a holomorphic function $h$ defined in a
neighborhood of $(f^1(m), \ldots, f^d(m)) \in \C^d$.  $\V$ is hypocomplex
on $M$ if it is hypocomplex at each $m \in M$.
\end{definition}
It is easily checked that this definition is independent of the choice of
independent smooth solutions $f^1, \ldots, f^d$.

Most involutive structures are not hypocomplex.  A second family of
hypocomplex examples consists of the CR structures
induced on hypersurfaces of $\C^d$
with Levi form having at least one positive and at least one negative
eigenvalue.  The coordinates $z^1, \ldots, z^d$ restrict to
form a set of independent
smooth solutions of $\V$, and by the H. Lewy extension theorem, any
solution of $\V$ (i.e.\ any CR distribution) is the restriction of a
holomorphic function defined in a neighborhood of the hypersurface.

We are interested in a particular involutive structure which arises
naturally on the real blow-up of $\R^n$ in $\C^n$.  Recall that if $N$ is a
smooth submanifold of another smooth manifold $M$, then there is
canonically associated another smooth manifold $B$, the real blow-up of $N$
in $M$, together with a blow-down map $b: B \to M$.  The space $B$ is
constructed by replacing each
$n \in N$ by the projectivized normal space $\bP(TM/TN)$ at $n$.  The
set $\Sigma \equiv b^{-1}(N)$ is a hypersurface in $B$.
The map $b:B \setminus \Sigma \to M
\setminus N$ is a diffeomorphism but, for $n \in N$, the fiber $b^{-1}(n)$
is diffeomorphic to
$\R\bP^{d-1}$, where $d$ is the codimension of $N$ in~$M$.

The prototypical
example is the blow-up of $\{0\} \subset \R^n$.  In this case local
coordinates $(y,t^1, \ldots, t^{n-1})$ on $B$ can be obtained in a
neighborhood of the point of the fiber above $\{0\}$ determined by the line
through the first basis vector by writing a point of $\R^n$ with nonzero first
coordinate in the form $(y,yt^1, \ldots, yt^{n-1})$.  If we set $t=(t^1,
\ldots, t^{n-1}) \in \R^{n-1}$, then the blow-down map $b$ takes the form
$$\begin{CD}
(y,t) & \quad\in\quad & B \\
\mapsdownto & & @VV{b}V \\
(y,yt)   & \quad\in\quad & \R^n
\end{CD}$$
and in these coordinates $\Sigma$ becomes $\{y=0\}$.
The image of this chart under $b$
consists of the origin together with all points with nonzero first
coordinate.  The full fiber
$\Sigma = b^{-1}(0)$ can be covered by charts obtained similarly using
each of the standard basis directions.

Coordinates for a general blow-up of $N \subset M$ can be obtained by
choosing appropriate local coordinates on $M$, applying the above
construction for the coordinates transverse to $N$, and leaving the
coordinates along $N$ unchanged.  For the case of interest here, that of
$\R^n \subset \C^n$, this amounts to
applying the above where $(y,t)$ are the imaginary parts of the variables.
So we write the real parts as $(x,s)$ for $x \in \R$ and $s \in \R^{n-1}$
and set $z=x+iy \in \C$.  It will turn out to be convenient to take our
blow-down map to be
a slight modification of the one just described, namely
$$\begin{CD}
(z,s,t) & \quad\in\quad & \C\times\R^{n-1}\times\R^{n-1} \subset B \\
\mapsdownto   \\
(z,s+zt)   & \quad\in\quad & \C\times\C^{n-1} \cong \C^n.
\end{CD}$$
This map is still exactly the standard blow-down in the imaginary
parts and for fixed imaginary parts is simply an invertible linear
transformation in
the real parts, so certainly provides local coordinates realizing
$b:B \to \C^n$.  In these coordinates $\Sigma$ again becomes $\{y=0\}$ and
$\Sigma$
is covered by charts obtained similarly using each of the standard basis
directions in $\R^n$.

Since $b:B \setminus \Sigma \to \C^n \setminus \R^n$ is a diffeomorphism,
the complex structure on $\C^n \setminus \R^n$ pulls back to a complex
structure on $B \setminus \Sigma$.  Although this complex structure does
not extend across $\Sigma$, it does extend as an involutive structure.
\begin{proposition}
The involutive structure on $B \setminus \Sigma$ obtained by pulling back
the complex structure on $\C^n \setminus \R^n$ extends smoothly across
$\Sigma$ to determine an involutive structure on all of $B$.
\end{proposition}
\begin{proof}
Denote points in $\C^n$ by $(z,w)$ with $z \in \C$ and $w \in \C^{n-1}$.
The complex
structure on $B \setminus \Sigma$ is given by $\V = \ker \{b^*dz, b^*dw\}$,
where we interpret $dw$ as vector-valued.  Using the coordinate expression
for $b$ given above, this is $\V = \ker \{dz, ds + t dz + z dt\}$.  As these
forms clearly remain linearly independent across $\Sigma$, $\V$ extends
smoothly across $\Sigma$ as a vector bundle, so it follows that
it determines an involutive structure on $B$.
\end{proof}

{From} the above it follows that in our $(z,s,t)$  coordinates on $B$, $\V$
is spanned by $\partial_{\zb}$ and the $\partial_t - z \partial_s$.
On $\Sigma$ the latter are real and span $\V \cap \overline{\V}$, which
agrees with the complexified tangent space to the fibers of $b|_{\Sigma}$.
In particular $\V \cap \overline{\V}$ has dimension $n-1$ on $\Sigma$.
The coordinates $z$ and $w=s+zt$ are solutions of $\V$, so
$\V$ is locally integrable.

For $n=2$, this involutive structure was introduced in \cite{hj}
in different coordinates as an example
as $V_2$, p.~501, and as (1.10), p.~502.

The main result of this section is the following.
\begin{theorem}\label{hypo}
The involutive structure on $B$ is hypocomplex.
\end{theorem}
\begin{proof}
Away from $\Sigma$ the involutive structure is complex, so certainly it is
hypocomplex there.  Hence it suffices to establish hypocomplexity at
$(z,s,t)=(0,0,0)$ in the coordinate system introduced above.  It must be
shown
that every distribution solution of $\V$ near $(z,s,t)=(0,0,0)$ is of the
form $h(z,s+zt)$ for some holomorphic function $h$ defined in a
neighborhood of  $(z,w)=(0,0) \in \C^n$.  We present the details for $C^1$
solutions of $\V$ and then indicate the modifications required to extend
the argument to distribution solutions.

Suppose that $f(z,s,t)$ is a $C^1$ solution of $\partial_{\zb}f=0$ and
$(\partial_t - z \partial_s)f=0$ in a neighborhood of the origin.  Then $f$
is, in particular, holomorphic in $z$, and so may be expanded as
\begin{equation}\label{start}
f(z,s,t)=\sum_{k=0}^{\infty} a_k(s,t)z^k,
\end{equation}
where the $a_k$ are $C^1$ functions near $(s,t) = (0,0)$, and this series
along with all first derivatives converges uniformly for $(z,s,t)$ small.
In particular, there are $C$, $r$ and $\epsilon > 0$ so that
\begin{equation}\label{cauchy}
|a_k(s,t)| \leq Cr^{-k} \qquad \text{for $|s|, |t| \leq \epsilon.$}
\end{equation}

Now $\partial_t f = z \partial_s f$ gives
$$\partial_t a_0 = 0 \qquad \text{and} \qquad \partial_t a_k = \partial_s
a_{k-1} \quad \text{for $k \geq 1$.}$$
{From} the first equation it follows that $a_0(s,t) = b_0(s)$ for a $C^1$
function $b_0(s)$.  {From} the second equation with $k=1$ one then obtains
$a_1(s,t)= t\cdot \partial_s b_0(s) + b_1(s)$ for a $C^1$
function $b_1(s)$.  This gives also that
$\partial_s b_0$ is $C^1$, so $b_0$ is $C^2$.  {From} the equations
with higher $k$ one similarly deduces inductively that
\begin{equation}\label{poly}
a_k(s,t) = \sum_{|\alpha|\leq
k}\frac{1}{\alpha!}t^{\alpha}\partial_s^{\alpha} b_{k-|\alpha|}(s)
\end{equation}
for $C^1$ functions $b_k(s)$, and the sum is over multiindices $\alpha$.
Since $a_k$ is $C^1$ it follows that each $\partial_s^{\alpha}
b_{k-|\alpha|}$ is $C^1$, so by induction each $b_k$ is actually
$C^{\infty}$.

In order to deduce bounds on the $b_k$ we use the following lemma, a
real version of Cauchy estimates for polynomials.
\begin{lemma}\label{kens}
For each $n$ there is a constant $R>0$ so that if $p(t) = \sum_{|\alpha|
\leq k} c_{\alpha}t^{\alpha}$ is a polynomial on $\R^n$ of degree $k$,
then
$$\text{max}\, |c_{\alpha}| \leq R^k \text{sup}_{|t_j| \leq 1}\,|p(t)|.$$
Here the sup is over all $t \in \R^n$ with $|t_1| \leq 1, \ldots, |t_n|
\leq 1$.
\end{lemma}
Lemma~\ref{kens} can be proved in a variety of ways.  The case
$n>1$ follows easily from the case $n=1$ by induction.  For $n=1$, the
coefficients of $p$ can be recovered from its values at $k+1$ points; taking
these points to be equally spaced and spread over $[-1,1]$ and making
straightforward estimates leads to a bound of the desired form (this
argument was
shown to us by Ken Bube).  Alternatively, for fixed $k$ and $l$,
$\sup|c_l|$ is explicitly
evaluated in terms of Chebyshev coefficients by a theorem of Markov, where
the sup is over all polynomials $p(t) = \sum_{l=0}^k c_l t^l$ satisfying
$\sup_{|t| \leq 1} |p(t)| \leq 1$ (see \cite{b}, p.~248).  A bound
as in Lemma~\ref{kens} follows easily from this.  We leave the details to
the reader.  By rescaling it is clear that in the statement of
Lemma~\ref{kens}, $\text{sup}_{|t_j| \leq 1}$ can be replaced by
$\text{sup}_{|t_j| \leq \epsilon}$, but now $R$ depends also on
$\epsilon > 0$.

By (\ref{poly}), $a_k$ is a polynomial in $t$ of degree $k$; applying
Lemma~\ref{kens} and using (\ref{cauchy}) one deduces the existence of
constants $C$ and $M$ so that for all $k$ and $\alpha$,
$$\sup_{|s| \leq \epsilon} \frac{1}{\alpha!}|\partial^{\alpha}b_k(s)| \leq
CM^{k+|\alpha|}.$$
In particular, each $b_k$ is a real-analytic function of $s$ and if we
expand
\begin{equation}\label{ra}
b_k(s) = \sum_{\alpha \geq 0} \frac{1}{\alpha!}c_{k,\alpha} s^{\alpha},
\end{equation}
then $\frac{1}{\alpha!} |c_{k,\alpha}| \leq CM^{k+|\alpha|}.$  Upon
substituting (\ref{ra}) into (\ref{poly}) and then into (\ref{start}) and
simplifying, one obtains
$$f=\sum_{\substack{k \geq 0\\ \alpha \geq 0}}
\frac{1}{\alpha!}c_{k,\alpha} z^k(s+zt)^{\alpha}.$$
Thus $f$ is of the desired form with
$h(z,w)=\sum \frac{1}{\alpha!}c_{k,\alpha} z^k w^{\alpha}.$

The same sort of argument applies for distribution solutions except that
the estimates must be made in different norms.  We sketch the steps and
leave the details to the reader. If $f$ is a distribution
solution then $f$ is still holomorphic in $z$, so can be expanded in the
form (\ref{start}), where now the $a_k$ are distributions in $(s,t)$ and the
series converges weakly in $(z,s,t)$.  Also in some neighborhood of the
origin $f$ has finite order, so (\ref{cauchy}) can be replaced by the
statement that for some $\epsilon$, $N$, and $r$, the set $\{r^k a_k\}$
defines a bounded set in $(C^N_{B_{\epsilon} \times B_{\epsilon}})'$.
Here $C^N_{B_{\epsilon} \times B_{\epsilon}}$ denotes the space of
restrictions to $B_{\epsilon} \times B_{\epsilon} = \{|s| \leq \epsilon\}
\times \{|t| \leq \epsilon\}$ of $C^N$ functions on $\R^{n-1} \times
\R^{n-1}$ with compact support in $B_{\epsilon} \times B_{\epsilon}$, with
the usual norm.  The inductive calculations and bootstrap leading to
(\ref{poly}) proceed just as
before, and it follows that the $a_k$ are of the form (\ref{poly}) for
$C^{\infty}$ functions $b_k(s)$.  It is not hard to prove a version of
Lemma~\ref{kens} for the $(C^N)'$ norm; the statement is the same except
that the right hand side is replaced by $C R^k \|p\|_{(C^N_{B_1})'}$, where
$C$ is another constant depending only on $n$.  Applying this version of
the Lemma now gives the bound
$$\frac{1}{\alpha!}\|\partial^{\alpha}b_k\|_{(C^N_{B_{\epsilon}})'} \leq
CM^{k+|\alpha|}.$$
{From} this it follows that each $b_k$ is real-analytic and that $f$ is of
the desired form for a holomorphic function $h(z,w)$ as before.
\end{proof}
\begin{remark}
F. Treves has pointed out to us that the hypocomplexity
also follows from microlocal regularity results for solutions of
involutive structures given in \cite{bt} or \cite{bct}.
The result in \cite{bt}
characterizes the analytic wave-front set of solutions of involutive
structures of tube type; the involutive structure on the blow-up is of this
form.  The result of \cite{bct} is that a point of the characteristic
variety of an involutive structure is not in the analytic wave-front set of
any solution if the Levi form at that point has a negative eigenvalue.  In
our case, at each point of the characteristic variety the Levi form has
exactly one positive and one negative eigenvalue.
\end{remark}

Next we show that the first local cohomology group of the involutive
structure is infinite dimensional at a point of $\Sigma$.  This amounts to
studying solvability of the inhomogeneous equations
$$\partial_{\zb} f = u\quad\mbox{and}\quad
(\partial_{t^j} - z \partial_{s^j})f = v_j,$$
subject to the compatibility
conditions
$$(\partial_{t^j} - z \partial_{s^j})u = \partial_{\zb} v_j\quad\mbox{and}
\quad(\partial_{t^j} - z \partial_{s^j})v_k
                   = (\partial_{t^k} - z \partial_{s^k})v_j.$$
We will take $u=0$ and $v_j=v_j(s)$ to depend only on~$s$.  In this case the
first compatibility condition is automatic and the
second reduces to $\partial_{s_j} v_k = \partial_{s_k} v_j$, which just
says that the 1-form $v_j ds^j$ is closed.  We assume {\em a priori} only
that the $v_j$ are distributions in $s$.

\begin{theorem}
If the equations $\partial_{\zb} f = 0$, $(\partial_{t^j} - z
\partial_{s^j})f = v_j(s)$ have a distribution solution $f$ near
$(z,s,t)= (0,0,0)$, then the $v_j(s)$ are real analytic near $s=0$.
\end{theorem}
\begin{proof}
The argument is similar to the proof of Theorem~\ref{hypo}.  Again we
present the details for the case in which there is a $C^1$ solution $f$;
the modifications to deal with a distribution solution follow
those outlined above.

The solution $f$ must be holomorphic in $z$, so we have
(\ref{start}) and (\ref{cauchy}) just as before.  The second equation gives
$$\partial_{t^j} a_0 = v_j(s) \qquad \text{and} \qquad \partial_{t^j} a_k
= \partial_{s^j} a_{k-1} \quad \text{for $k \geq 1$.}$$
It follows that $a_0(s,t) = t^j v_j(s) + b_0(s)$ for a function $b_0(s)$;
since $f$ is $C^1$ we deduce that $v_j$ and $b_0$ are both $C^1$ as well.
Inductively solving the remaining equations gives
\begin{equation}\label{poly2}
a_k(s,t) = \sum_{|\alpha|=k+1} \frac{1}{\alpha!}t^{\alpha}
(\partial_s^k v)_{\alpha} +
\sum_{|\alpha|\leq
k}\frac{1}{\alpha!}t^{\alpha}\partial_s^{\alpha} b_{k-|\alpha|}
\end{equation}
for $C^1$ functions $b_k(s)$,
where $(\partial_s^k v)_{\alpha}$ denotes the $\alpha$-multiindex
component of the symmetric tensor $\partial_s^k v$.  Again, one deduces
inductively that $v$ and all the $b_k$ are $C^{\infty}$.  Now apply
Lemma~\ref{kens} to the highest order coefficients of $a_k$ as a polynomial
in $t$ and use
(\ref{cauchy}) to deduce that for some constants $C$ and $M$,
$\sup_{|s|\leq \epsilon} \frac{1}{\alpha!}
|\partial_s^{\alpha}v_j(s)| \leq CM^{|\alpha|}$.  It follows that
each $v_j$ is real-analytic.
\end{proof}

\section{The Edge of the Wedge Theorem}\label{eow}
The hypocomplexity of the involutive structure on the blow-up of $\R^n$ in
$\C^n$ is, for all intents and purposes, equivalent to the Edge of the
Wedge Theorem.  (See \cite{vzs} for a discussion of the history and
development of this result.)

\begin{EOWT}\label{eowt}
Let $E \subset \R^n$ be open and let $C \subset \R^n$ be an open cone such
that $C = -C$.
Set $W=E+iC \subset \C^n$.  Let $\N$ be
an open set in $\C^n$ with $\N \cap \R^n = E$ and let $f$ be a
holomorphic function on $\N \cap W$.
Suppose that $\lim_{C \ni Y \rightarrow 0} f(\cdot + iY)$ exists
in $\D'(E)$, and is independent of how $Y \rightarrow 0$.  Then there is
another open set $\N' \subset \N$ with $\N' \cap \R^n = E$ and
a holomorphic function $h$ on $\N'$ such that $h=f$ on $\N' \cap W$.
\end{EOWT}

The geometry is simplest in case $n=1$.  Setting $\N_+ = \N \cap
\{y>0\}$ and $\N_- = \N \cap \{y<0\}$, one has holomorphic functions $f_+$
on $\N_+$ and $f_-$ on $\N_-$ whose distribution boundary values agree on
$E$.  One proof of the Edge of the Wedge Theorem in this case
is to observe that $f_+$ and $f_-$ together
define a distribution on $\N$ which is a distribution solution of the
Cauchy-Riemann equation.  This observation follows from the following
(standard) formulation of distribution boundary values of holomorphic
functions:
\begin{equation}\label{bv}
\begin{gathered}
\parbox{4.4in}
{Let $E \subset \R$ and  $\N \subset \C$ be open with $\N \cap \R = E$,
and set $\N_+ = \N \cap \{y>0\}$.   Suppose that $f$ is a holomorphic
function on
$\N_+$ and one of the following limits exists:
\begin{enumerate}
\renewcommand{\theenumi}{(\roman{enumi})}
\item
$\tilde f = \lim_{\epsilon \rightarrow 0} f \,\,
\chi^{\vphantom{|}}_{\{y > \epsilon\}}$ in $\D'(\N)$.
\item
$f_0 = \lim_{\epsilon \rightarrow 0} f(\cdot + i\epsilon)$ in $\D'(E)$.
\end{enumerate}
Then the other exists, and $\partial_{\zb} \tilde f = \frac{i}{2} f_0(x)
\delta_0(y)$.
}
\end{gathered}
\end{equation}
For boundary values from $\N_-$ the corresponding equation for
$\partial_{\zb} \tilde f$
has a minus sign.  Hence $\partial_{\zb} (\tilde f_+ + \tilde f_-) = 0$,
and $h=\tilde f_+ + \tilde f_-$ is the required holomorphic function.

For completeness, we briefly recall the proof of (\ref{bv}).  If (i) exists,
then certainly also $\lim_{\epsilon \rightarrow 0} \partial_{\zb}
(f \,\, \chi^{\vphantom{|}}_{\{y > \epsilon\}})$ exists in $\D'(\N)$.  However,
$\partial_{\zb} (f \,\, \chi^{\vphantom{|}}_{\{y > \epsilon\}}) =
\frac{i}{2} f \,\delta_{\epsilon}(y)$, so it follows that (ii) exists, and
in the limit one obtains $\partial_{\zb} \tilde f = \frac{i}{2} f_0(x)
\delta_0(y)$.  On the other hand, if (ii) exists, then $y \rightarrow
f(\cdot + iy)$ defines a weakly continuous function of $y$ into $\D'(E)$ up to
$y=0$.  Using the uniform boundedness principle one deduces easily that
(i) exists.

Although the geometry of the wedge is more complicated in higher
dimensions, introduction of the blow-up reduces the geometry to exactly
that in dimension one; namely solutions (of the involutive structure) on
the two sides of a hypersurface with equal boundary values, and the above
proof then goes right through.  The main
observation is simply that the image under the blow-down map $b$ of a ball in
$B$ about a point of $\Sigma$ is essentially a localized wedge as occurs in
the Edge of the Wedge Theorem.
\begin{proof}
Choose a point of $E$ and a line in $C$; these determine a
line in $W$ and a point in the
fiber in $B$ above the point of $E$.  Coordinates may be chosen so that the
point is the origin in $\R^n$ and the line is the first basis direction,
so that in the coordinates introduced above the corresponding point of
$B$ is the origin $(z,s,t) = (0,0,0)$.
Let $\M$ be an open ball about the origin in $B$ and set $\M_+ = \M \cap
\{y>0\}$, $\M_- = \M \cap \{y<0\}$, $\M_0 = \M \cap \{y=0\}$.
Then for $\M$ sufficiently small, $b(\M_+ \cup \M_-) \subset \N \cap W$.
Consequently $f\circ b$ defines solutions $f_+$
and $f_-$ of the involutive structure on $\M_+$ and $\M_-$, resp.  In the
case that the boundary values of $f$ are taken continuously, so that $f$
has a continuous extension to $E \cup (\N \cap W)$, $f\circ b$ clearly
extends continuously to $\M$ so the boundary values of $f_+$ and $f_-$
agree on $\M_0$.  In the general case in which the boundary values of $f$
are taken in the distribution sense, it is straightforward to check (again
using the uniform boundedness principle) that
the distribution boundary values of $f_+$ and $f_-$ exist in $\D'(\M_0)$ and
agree.  Since $f_+$ and $f_-$ are
annihilated by $\partial_{\zb}$, one can now apply the exact argument
as presented above in the $n=1$ case to conclude that
$\partial_{\zb}(\tilde f_+ + \tilde f_-) = 0$.  However
$(\partial_t - z \partial_s)\tilde f_+=
\lim (\partial_t - z \partial_s)(f_+ \,\, \chi^{\vphantom{|}}_{\{y >
\epsilon\}}) = 0$, and similarly for $\tilde f_-$.  Thus
$\tilde f_+ + \tilde f_-$ is a distribution solution
of the involutive structure in $\M$.  {From} Theorem~\ref{hypo} it follows
that there is a holomorphic function in a neighborhood of $0 \in \C^n$
which agrees with $f$ on $b(\M') \setminus E$
for some ball $\M' \subset \M$.  In
particular, $f$ extends continuously up to $\R^n$.
Since a holomorphic function is determined by its
restriction to $\R^n$ and the boundary values are independent of the
direction chosen, it follows that the holomorphic extensions
obtained by varying the point of $E$ and the direction in $C$ all agree on
overlaps, so they define one holomorphic function in a
neighborhood of $E$ as desired.
\end{proof}

It is also easy to deduce the hypocomplexity of the involutive structure
from the Edge of the Wedge Theorem.  Since on $\Sigma$, the $\partial_t - z
\partial_s$ are real and span the tangent space to the fibers of
$b|_{\Sigma}$, any $C^1$ solution of the involutive structure on a ball
about a point of $\Sigma$ is constant on these fibers so defines a
holomorphic function on a localized wedge, continuous up to $\R^n$.  The
Edge of the Wedge Theorem provides the holomorphic extension
to a neighborhood as required in
the definition of hypocomplexity.  This same argument works for
distribution solutions as well, using the fact that any distribution
solution of the involutive structure is transversally regular, i.e.\ it
can be regarded as a continuous
function of $y$ with values in distributions on a piece of $\Sigma$.

\section{Integral Geometry}\label{intgeom}
Consideration of the involutive structure on the blow-up of $\R\bP^n$ in
$\C\bP^n$ arose in \cite{begm}, \cite{e}, and \cite{be} in applying methods of
complex integral geometry, specifically the Penrose transform,
to problems in real integral geometry.  We include here a
description of the underlying geometric picture.  This provides a natural
realization of the blow-up $B$ in this compactified setting.

The real problem is to study the X-ray transform in
$\R^n$, or really a compactified version thereof.  Embed
$\R^n \subset \R\bP^n$; then an affine line in $\R^n$ is represented as a
2-plane
through the origin in $\R^{n+1}$, so the space of such lines becomes a
subset of $\mathrm{Gr}_2(\R^{n+1})$.
The correspondence is encoded in the double fibration
$$\begin{picture}(250,60)
\put(58,50){\makebox(0,0){$\F_{1,2}(\R^{n+1})$}}
\put(40,40){\vector(-3,-4){15}}
\put(60,40){\vector(3,-4){15}}
\put(20,10){\makebox(0,0){$\R\bP^n$}}
\put(95,10){\makebox(0,0){$\mathrm{Gr}_2(\R^{n+1})$,}}
\end{picture}$$
where $\F_{1,2}(\R^{n+1})$ denotes the flag manifold of lines contained in
planes in $\R^{n+1}$.  The analogous complex correspondence is
$$\begin{picture}(250,60)
\put(58,50){\makebox(0,0){$\F_{1,2}(\C^{n+1})$}}
\put(40,40){\vector(-3,-4){15}}
\put(60,40){\vector(3,-4){15}}
\put(20,10){\makebox(0,0){$\C\bP^n$}}
\put(95,10){\makebox(0,0){$\mathrm{Gr}_2(\C^{n+1})$.}}
\put(27,32){\makebox(0,0){$\mu$}}
\put(73,32){\makebox(0,0){$\nu$}}
\end{picture}$$
The complex methods are
applied to the real problem via the introduction of a hybrid
correspondence.  Define $F=\nu^{-1}(\mathrm{Gr}_2(\R^{n+1}))$, so that
$$F=\{(L,P): L\subset \C^{n+1}\,\, \text{is a complex line,}\,\, P\subset
\R^{n+1}\,\, \text{is a real 2-plane,}$$
$\qquad \qquad \qquad \,\, \text{and}\,\, L \subset P+iP\}.$
\vspace{.1 in}

\noindent
The hybrid correspondence is then
$$\begin{picture}(250,60)
\put(50,50){\makebox(0,0){$F$}}
\put(40,40){\vector(-3,-4){15}}
\put(60,40){\vector(3,-4){15}}
\put(3,10){\makebox(0,0){$\R\bP^n \subset \C\bP^n$}}
\put(93,10){\makebox(0,0){$\mathrm{Gr}_2(\R^{n+1})$.}}
\put(22,36){\makebox(0,0){$\mu|_F$}}
\put(80,36){\makebox(0,0){$\nu|_F$}}
\end{picture}$$

Of course $\nu|_F$ is still a fibration, but this is not the case for
$\mu|_F$.  Write $L \in \C\bP^n$ as $L=\,\, \text{span}_{\C}\,\{\xi +
i\eta\}$ for $\xi, \eta \in \R^{n+1}$.  If $L \notin \R\bP^n$ then $\xi$
and $\eta$ are linearly independent, so if $L \subset P+iP$ for a real
2-plane $P$, necessarily $P=\,\,\text{span}\{\xi,\eta\}$.  Hence
$(\mu|_F)^{-1}(L)$ consists of a single point.  On the other hand, if
$L \in \R\bP^n$, then $L=\,\, \text{span}_{\C}\,\{\xi\}$ for some
$\xi \in \R^{n+1}$ so $\mu^{-1}(\R\bP^n)=\F_{1,2}(\R^{n+1})$ and
$(\mu|_F)^{-1}(L) \cong \R\bP^{n-1}$.
This suggests the following:

\begin{proposition}
The mapping $\mu|_F : F \rightarrow \C\bP^n$ realizes the blow-up of
$\R\bP^n \subset \C\bP^n$.
\end{proposition}
\begin{proof}
We show that, in suitable coordinates, $\mu|_F$ is
exactly the map $b$ introduced previously.  Coordinates in $\R^{n+1}$ may
be chosen to represent any fixed point of $\F_{1,2}(\R^{n+1})$ as
$(\text{span}\{e_1\},
\text{span}\{e_1, e_2\})$, where $e_1, e_2$ are the first two basis vectors.
For an arbitrary nearby point \mbox{$(L,P) \in F$}, we have
$$P=\,\text{span}_{\R}
\left\{
\left[
\begin{matrix}
1\\
0\\
s
\end{matrix}
\right],
\left[
\begin{matrix}
0\\
1\\
t
\end{matrix}
\right]
\right\}\quad
\mbox{for some $s,t \in \R^{n-1}$,}$$
and
$$L=\,\,\text{span}_{\C}
\left\{
\left[
\begin{matrix}
1\\
0\\
s
\end{matrix}
\right]
+z
\left[
\begin{matrix}
0\\
1\\
t
\end{matrix}
\right]=
\left[
\begin{matrix}
1\\
z\\
s+zt
\end{matrix}
\right]
\right\}\quad
\mbox{for some $z \in \C$}.$$
Coordinates on $F$ are $(z,s,t)$, and $\mu|_F$ becomes
$(z,s,t) \rightarrow (z,s+zt)$ as desired.
\end{proof}

We refer to \cite{begm}, \cite{e}, and \cite{be} for the application of this
geometry to the inversion of real integral transforms.


\begin{thebibliography}{BEGM}

\bibitem[BaE]{be} T.N. Bailey and M.G. Eastwood,
{\em Zero-energy fields on real projective space},
Geom. Dedicata, to appear.

\bibitem[BEGM]{begm} T.N. Bailey, M.G. Eastwood, A.R. Gover, and L.J. Mason,
{\em The Funk transform as a Penrose transform}, preprint.

\bibitem[BT]{bt}
M.S. Baouendi and F. Treves,
{\em A microlocal version of Bochner's tube theorem},
Indiana Math. Jour. {\bf 31} (1982), 885---895.

\bibitem[BCT]{bct}
M.S. Baouendi, C.H. Chang, and F. Treves,
{\em Microlocal hypo-analyticity and extension of CR functions},
Jour. Diff. Geom. {\bf 18} (1983), 331--391.

\bibitem[BoE]{b} P. Borwein and T. Erd\' elyi,
{\em Polynomials and Polynomial Inequalities},
Grad. Text. Math. vol. 161, Springer, 1995.

\bibitem[CT]{ct} P.D. Cordaro and F. Treves,
{\em Hyperfunctions on Hypo-analytic Manifolds},
Ann. Math. Stud. vol. 136, Princeton University Press, 1994.

\bibitem[E]{e} M.G. Eastwood,
{\em Complex methods in real integral geometry},
Proceedings of the Fifteenth Winter School on Geometry and Physics,
Srn\'\i\/ 1995,
Suppl. Rend. Circ. Mat. Palermo, to appear.

\bibitem[HJ]{hj} N. Hanges and H. Jacobowitz,
{\em Involutive structures on compact manifolds},
Amer. Jour. Math. {\bf 177} (1995), 491--522.

\bibitem[T]{t} F. Treves,
{\em Hypo-analytic structures},
Princeton University Press, 1992.

\bibitem[VZS]{vzs}
V.S. Vladimirov, V.V. Zharinov, and A.G. Sergeev,
{\em Bogolyubov's ``edge of the wedge'' theorem, its development and
applications},
Russian Math. Surveys {\bf 49:5} (1994), 51--65.

\end{thebibliography}
\end{document}